\documentclass{amsart}
\usepackage{amsmath}
\usepackage{amsthm}
\usepackage{amssymb}
\usepackage{color}

\newtheorem{thm}{Theorem}
\newtheorem{lem}[thm]{Lemma}

\newtheorem{defi}[thm]{Definition}

\newtheorem{rk}[thm]{Remark}

\newcommand{\rr}{{\mathbb{R}}}
\newcommand{\rd}{{\rr^3}}
\newcommand{\nn}{{\mathbb{N}}}
\newcommand{\sS}{{\mathbb{S}^2}}
\newcommand{\e}{\epsilon}
\newcommand{\vip}{\vskip.2cm}
\newcommand{\indiq}{\hbox{\rm 1}{\hskip -2.8 pt}\hbox{\rm I}}
\newcommand{\intrd}{\int_{\rr^3}}
\newcommand{\dd}{{\rm d}}
\newcommand{\vs}{{v_*}}
\newcommand{\cA}{{\mathcal{A}}}
\newcommand{\cB}{{\mathcal{B}}}
\newcommand{\cC}{{\mathcal{C}}}
\newcommand{\cP}{{\mathcal{P}}}
\newcommand{\ig}{[\![}
\newcommand{\id}{]\!]}

\newcommand{\bla}{\color{black}}

\begin{document}

\title[Homogeneous Boltzmann equations]
{On exponential moments of the 
homogeneous Boltzmann equation for hard potentials without cutoff}

\author{Nicolas Fournier}
\address{Sorbonne Universit\'e, LPSM-UMR 8001, Case courrier 158,75252 Paris Cedex 05, France.}
\email{nicolas.fournier@sorbonne-universite.fr}


\begin{abstract} 
We consider the spatially homogeneous Boltzmann equation for hard potentials without cutoff.
We prove that an exponential moment of order $\rho=\min\{2\gamma/(2-\nu),2\}$, 
with the usual notation, is immediately created. This is stronger than what happens in the case
with cutoff.
We also show that exponential moments of order $\rho\in (0,2]$ are propagated.
\end{abstract}

\maketitle

\section{Introduction and results}
\setcounter{equation}{0}

\subsection{The Boltzmann equation}

We consider a spatially homogeneous gas modeled by the Boltzmann equation:
the density $f_t(v)$ of particles with velocity $v\in \rr^3$ at time $t\geq 0$ solves
\begin{gather} \label{be}
\partial_t f_t(v) = \intrd \dd \vs  \int_{\sS} \dd\sigma B(|v-\vs |,\cos \theta)
\big[f_t(v')f_t(v'_*) -f_t(v)f_t(\vs )\big],
\end{gather}
where
\begin{gather}
\label{vprime}
v'=\frac{v+\vs }{2} + \frac{|v-\vs |}{2}\sigma, \quad 
v'_*=\frac{v+\vs }{2} -\frac{|v-\vs |}{2}\sigma \quad \hbox{and} 
\quad \cos \theta=  \frac{v-\vs }{|v-\vs |} \cdot \sigma  .
\end{gather}
We refer to the book of Cercignani \cite{c} and to the long review papers of Villani \cite{vh} and 
Alexandre \cite{a} for some detailed and complete accounts of what is known, both from the physical
and mathematical points of view, about this equation.
One may assume without loss of generality that the initial condition satisfies
$$
\intrd f_0(v)\dd v = 1, \quad \intrd vf_0(v)\dd v = 0 \quad \hbox{and}\quad \intrd |v|^2 f_0(v)\dd v = 1,
$$
and these quantities, namely the mass, momentum and kinetic energy, 
are constant, at least informally, as time evolves.

\subsection{Assumptions}
We will suppose that for some $\gamma\in (0,1]$ and some $\nu\in(0,2)$,
\begin{align*}
&(H_1(\gamma)) \qquad B(|v-\vs |,\cos \theta)\sin \theta = |v-\vs |^\gamma \beta(\theta)
\quad \hbox{for some measurable $\beta:(0,\pi]\to\rr_+$,}\\[2mm]
&(H_2(\nu)) \qquad\exists \; \kappa_1,\kappa_2\in (0,\infty),\;  
\forall \; \theta \in (0,\pi], \;\; \kappa_1 \theta^{-\nu-1} \leq \beta(\theta) \leq \kappa_2  \theta^{-\nu-1} .
\end{align*}

As explained in \cite{c,vh,a}, when particles interact through a repulsive force 
in $1/r^s$ for some $s> 2$, we have $(H_1(\gamma))$ and $(H_2(\nu))$ with 
$\gamma=(s-5)/(s-1)$ and $\nu=2/(s-1)$.
When $\gamma \in (0,1)$ (i.e. $s>5$), one speaks of hard potentials.

\vip

One speaks of hard potentials with cutoff when we have $(H_1(\gamma))$ for some $\gamma \in (0,1]$
and when $(H_2(\nu))$ is replaced by the condition
$\int_0^\pi \beta(\theta)\dd\theta \in (0,\infty)$, which more or less corresponds
to the case where $\nu=0$.

\subsection{Weak solutions}\label{wsp}

First, we parameterize (\ref{vprime}) as in \cite{fme}.
For every $X\in \rr^3\setminus \{0\}$, we introduce $I(X),J(X)\in\rr^3$ such that
$(\frac{X}{|X|},\frac{I(X)}{|X|},\frac{J(X)}{|X|})$ 
is an orthonormal basis of $\rr^3$. We also put $I(0)=J(0)=0$.
For $X,v,\vs \in \rr^3$, $\theta \in (0,\pi]$ and $\varphi\in[0,2\pi)$,
we set
\begin{align}\label{rep}
\left\{
\begin{array}{l}
\Gamma(X,\varphi)=(\cos \varphi) I(X) + (\sin \varphi)J(X), \\[2mm]
\displaystyle v'= v - \frac{1-\cos\theta}{2} (v-\vs )
+ \frac{\sin\theta}{2}\Gamma(v-\vs ,\varphi),\\[2mm]
\displaystyle v'_*= \vs + \frac{1-\cos\theta}{2} (v-\vs )
- \frac{\sin\theta}{2}\Gamma(v-\vs ,\varphi).
\end{array}
\right.
\end{align}
We denote by $\cP(\rd)$ the set of probability measures on $\rd$.
For $p\in \rr_+$ and $f\in \cP(\rd)$, we introduce the moment of order $p$ of $f$: 
$$m_p(f)=\intrd |v|^p f(\dd v).$$
We use the following classical notion of weak solutions.

\begin{defi}\label{dfws}
Assume $(H_1(\gamma))$ and $(H_2(\nu))$ for some $\gamma \in (0,1]$ and some $\nu\in(0,2)$.
A weakly continuous family $(f_t)_{t\geq 0}$ of probability measures on $\rd$
is a weak solution to (\ref{be}) if for all $t\geq 0$,
\begin{align*}
\intrd vf_t(\dd v)=0 \quad \hbox{and} \quad m_2(f_t)=1
\end{align*}
and if for any $\phi \in C^2_b(\rd)$ and any $t\geq 0$, using the parameterization \eqref{rep},
\begin{align}\label{weak}
\frac{\dd}{\dd t}\!\intrd \!\!\phi(v)f_t(\dd v)\! =\!\!
\intrd \!\intrd\! \int_0^\pi \!\!\!\int_0^{2\pi} \! [\phi(v')\!+\phi(\vs')\!-\phi(v)\!-\phi(\vs)]|v\!-\!\vs|^\gamma
\dd \varphi \beta(\theta)\dd \theta f_t(\dd \vs) f_t(\dd v).
\end{align}
\end{defi}

As shown by Lu-Mouhot in \cite{lm}, see also Villani \cite{vh}, weak solutions exist
starting from any given initial condition $f_0\in \cP(\rd)$ such that $\intrd v f_0(\dd v)=0$
and $m_2(f_0)=1$, and they satisfy
\begin{gather}\label{momnul1}
 \textstyle\hbox{for all $p\geq 0$, all $t_0>0$,} \quad \sup_{t\geq t_0} m_p(f_t)<\infty.
\end{gather}
Let us mention the recent uniqueness result of Heydecker \cite{h}, which concerns the
case where $\nu\in(0,1)$,
assuming only that $m_p(f_0)<\infty$  for some sufficiently large $p$.

\subsection{Main result}

Here is our main result.

\begin{thm}\label{mr}
Assume $(H_1(\gamma))$ and $(H_2(\nu))$ for some $\gamma \in (0,1]$ and some
$\nu \in (0,2)$. 
Consider any weak solution $(f_t)_{t\geq 0}$ to \eqref{be}.
\vip
(i) Put $\rho=\min\{2\gamma/(2-\nu),2\}$. There are some constants $T>0$ and $\sigma>0$,
depending only on $\gamma,\nu,\kappa_1,\kappa_2$, such that 
$$
\sup_{t\in [0,T]} \intrd \exp[\sigma t^{\rho/\gamma} |v|^{\rho}] f_t(\dd v) \leq 4.
$$

(ii) For any $A>0$, any $\sigma_0>0$, any $\rho\in (0,2]$, there is a constant
$\sigma>0$, depending only on $\gamma,\nu,\kappa_1,\kappa_2,\rho,\sigma_0,A$, such that
$$
\intrd \exp[\sigma_0 |v|^{\rho}] f_0(\dd v)\leq A \quad \Longrightarrow \quad
\sup_{t\geq 0} \intrd \exp[\sigma |v|^{\rho}] f_t(\dd v) \leq 6.
$$
\end{thm}

Since $\min\{2\gamma/(2-\nu),2\}>\gamma$, point (i) is stronger than in 
the cutoff case where, as we will see in the next subsection, only exponential moments of order $\rho=\gamma$
are created.
\vip
In (ii), we have a possible deterioration of 
the constant $\sigma$, as in all the references below.

\vip
By (ii), (i) can be extended to: there is $\sigma>0$ such that, with $\rho=\min\{2\gamma/(2-\nu),2\}$, 
$$
\sup_{t\geq 0} \intrd \exp[\sigma \min\{t^{\rho/\gamma},1\} |v|^{\rho}] f_t(\dd v) \leq 6.
$$

\subsection{References}
There is a large literature on the subject, because exponential moments can be used for different
purposes, such as estimating the rate of convergence to equilibrium, see Mouhot \cite{m}, or uniqueness,
see \cite{fm}.

\vip

Using the famous Povzner inequality \cite{p}, Wennberg \cite{w} discovered that polynomial moments
are immediately created by the homogeneous Boltzmann equation for cutoff hard potentials (CHP in short),
i.e. when the angular cross section $\beta$ is assumed to be
integrable on $[0,\pi]$, which roughly corresponds to the case $\nu=0$.
This really requires that $\gamma>0$ and the main intuition is that particles with large
velocities are quickly slowed down, because they collide at large rate 
(since $\gamma>0$) with slow particles.

\vip

In his seminal paper \cite{b}, Bobylev proved that Gaussian moments ($\rho=2$) are propagated, still for CHP.
This relies on very tight computations involving a recursive ODE argument to estimate
the polynomial moments, which are then summed to estimate Gaussian moments.
Let us also cite Bobylev-Gamba-Panferov \cite{bgp} who studied inelastic collisions.
Following the ideas of \cite{b}, Mouhot \cite{m} managed to create exponential moments of order
$\rho=\gamma/2$, still for CHP. Following the same approach, Lu-Mouhot \cite{lm}
were able to create exponential moments of order $\rho=\gamma$ for CHP and non-cutoff hard potentials (NCHP).
Alonso-Ca\~nizo-Gamba-Mouhot \cite{acgm} found a much simpler method
to create exponential moments of order $\rho=\gamma$ and propagate exponential
moments of order $\rho\in(0,2]$, for CHP.
Let us finally quote Alonso-Gamba-Taskovi\'c \cite{agt}, who studied some much stronger
Lebesgue and Sobolev norms
of $f_t$ with exponential weights for CHP.

\vip

Concerning NCHP, there is the work of Lu-Mouhot \cite{lm} already mentioned.
The proof of Lemma 4.1 in Fournier-Mouhot \cite{fm}, which concerns exponential moments for NCHP, 
is unfortunately false, there is a major gap
(the function $\delta(\eta)$ in (4.6) actually depends on $p$).
What was required there for the uniqueness criterion for NCHP to imply a well-posedness result, 
was the propagation of
exponential moments of order $\rho=\gamma$.
Taskovi\'c-Alonso-Gamba-Pavlovi\'c \cite{tagp} have shown, for NCHP, creation of
exponential moments of order $\rho=\gamma$ and propagation of exponential moments of order $\rho\in(0,4/(2+\nu)]$
(which contains $\rho=\gamma\in (0,1]$ and thus fixes the issue in \cite{fm}).

\vip

It thus seems that concerning NCHP, \cite{tagp} contains the best available results,
and Theorem \ref{mr} is stronger both for creation and propagation.
In particular, we show that NCHP create more exponential moments than CHP.

\vip

The homogeneous Landau equation for hard potentials, which often behaves in a similar way
as the Boltzmann equation, but 
which is considerably simpler in many points, immediately creates Gaussian moments ($\rho=2$) for
any value of $\gamma\in(0,1]$, see \cite{fh}. Once this is observed, it is natural to wonder
if NCHP create more exponential moments than CHP. 
The answer is not intuitively clear, because the (considerably many) additional collisions
caused by the singularity of $\beta$ near $0$
involve some (considerably) small values of $\theta$ and do not much slow down particles
with high velocity. Actually, the effect is strong enough to modify the behavior of the solutions:
exponential moments of order $\rho=\min\{2\gamma/(2-\nu),2\}>\gamma$ are created by NCHP,
while only exponential moments of order $\rho=\gamma$ are created by CHP.
\vip

By the way, Theorem \ref{mr}-(i) implies that
for any $\gamma\in (0,1]$, Gaussian moments ($\rho=2$) are created by NCHP if the angular cross section
is singular enough, namely if $\nu \geq 2-\gamma$.

\vip

Our strategy is the same as that of \cite{tagp}: we adapt the ideas of
\cite{acgm} to NCHP, taking advantage of the simplicity of 
the method. The present paper resembles \cite{tagp} in several points. In particular, 
some {\it Mittag-Leffler moments} of the form 
$$\sum_{n\geq 0}\frac{a^n m_{2n}(f)}{\Gamma(\alpha n+1)},$$
with $a>0$ and $\alpha \geq 1$,
are used in \cite{tagp}, while 
we are led to use some series of the form 
$$\sum_{n\geq 0} \frac{a^n m_{2n}(f)}{(n!)^\alpha},$$
with $a>0$ and $\alpha\geq 1$. This is almost the same thing and in both cases, 
this is enough to control some exponential
moments of the form $\intrd \exp(b |v|^\rho) f(\dd v)$,
with $\rho=2/\alpha$, see Lemma \ref{biznorm}.
The main advantage of using such series is that it then suffices to study integer moments
$m_{2n}(f)$, which leads to much more explicit
computations than if using non-integer moments, as is done e.g. in \cite{acgm},
where $\intrd \exp(b |v|^\rho) f(\dd v)$ is more naturally studied through 
$\sum_{n\geq 0} (n!)^{-1}b^n m_{\rho n}(f)$.
\vip

However, we try to really take advantage of the singularity of the cross section
to establish a stronger Povzner inequality than in the cutoff case, see Lemma \ref{pov}
and the paragraph below. We then have to adapt suitably the
proof of \cite{acgm}, on the one hand because we can only deal with
integer moments, and on the other hand because we have to exploit the new Povzner
inequality.

\subsection{About optimality}
Since the solutions to \eqref{be} converge to some Maxwell (Gaussian) distributions, which are 
stationary solutions, we cannot expect to
create or propagate exponential moments of order $\rho>2$. 
The propagation result thus seems optimal.
Concerning creation, one may get convinced, following the proofs of Lemmas \ref{pov} and \ref{2groc},
that for some constant $c>0$,
$$
\forall\, n\geq 2,\quad m_{2n}'(f_t) \geq -cn^{\nu/2}m_{2n+\gamma}(f_t).
$$
Admitting, {\it and this is not so clear}, that 
the H\"older inequality is sharp enough so that we have $m_{2n+\gamma}(f_t)\simeq [m_{2n}(f_t)]^{1+\gamma/2n}$, 
we end with
$m_{2n}'(f_t) \gtrsim - n^{\nu/2} [m_{2n}(f_t)]^{1+\gamma/2n}$,
from which we easily conclude, if $m_{2n}(0)=\infty$, that $m_{2n}(f_t) \gtrsim [n^{1-\nu/2}/ t]^{2n/\gamma}$.
Still informally, this should tell us that $m_{\rho n}(f_t) \gtrsim [n^{1-\nu/2}/ t]^{\rho n/\gamma}$, so that
$$
\intrd \exp[\sigma |v|^\rho]f_t(\dd v) = \sum_{n\geq 0} \frac{\sigma^n m_{\rho n}(f_t) }
{n!} \gtrsim \sum_{n\geq 0} \frac{\sigma^n n^{\rho(1-\nu/2)n/\gamma}}{t^{\rho n/\gamma} n!}.
$$
By Stirling's formula $n!\sim \sqrt{2\pi n}(n/e)^n$, 
this series is divergent, for any value of $\sigma>0$, when $\rho(1-\nu/2)/\gamma>1$, i.e. when
$\rho > 2\gamma/(2-\nu)$, which is coherent with Theorem \ref{mr}-(i).

\subsection{About uniqueness}
Assume $(H_1(\gamma))$ and $(H_2(\nu))$ for some $\gamma \in (0,1]$ and some $\nu \in (0,1)$. 
Using Theorem \ref{mr} and \cite[Theorem 2.2]{fm}, it seems possible to prove, in a few pages, 
the well-posedness
of \eqref{be} assuming that the
initial condition satisfies $\intrd \exp{(|v|^\delta)} f_0(\dd v)<\infty$ for some $\delta>0$.
This is stronger than \cite{fm}, where we assumed that  $\intrd \exp{(|v|^\gamma)} f_0(\dd v)<\infty$,
but weaker than the recent result of Heydecker \cite{h}, who only assumes that $m_p(f_0)<\infty$
for some large explicitable $p$.

\subsection{Plan}
The paper is technical and we are guided by computations rather than intuition.

\vip
In Section \ref{pv}, we establish a Povzner lemma, which is stronger than what is known
in the cutoff case. We handle the whole computation as explicitly as possible, not relying on
any previous Povzner estimate, because this is required if we 
really want to show that the singular part of the cross section
accelerates the slowing down of particles.

\vip

In Section \ref{smom}, we derive some differential inequalities for the even integer moment
from the Povzner inequality, and we prove some first rough estimates about these moments.

\vip

In Section \ref{sr}, we quickly study how to control exponential moments by
even integer moments and {\it vice-versa}.

\vip

Finally, we adapt the proofs of \cite{acgm} 
to show Theorem \ref{mr}-(i) in Section \ref{sc}
and Theorem \ref{mr}-(ii) in Section \ref{sp}. This requires some work,
because we can only use integer moments, and because we start from a different Povzner
estimate.

\subsection{Notation}

We use the convention that $\nn=\{0,1,...\}$. For $a,b\in \nn$ with $a\leq b$, we set
$\ig a,b \id=\{a,a+1,\dots,b\}$. In the whole paper, $(f_t)_{t\geq 0}$ is a given 
weak solution satisfying $m_0(f_t)=m_2(f_t)=1$ and $\intrd v f_t(\dd v)=0$ for all $t\geq 0$.
For $p\in\rr_+$ and $t\geq 0$, we set $m_p(t)=m_p(f_t)$.

\section{A non-cutoff Povzner lemma}\label{pv}

The goal of this section is to establish the following Povzner inequality.

\begin{lem}\label{pov}
Assume $(H_2(\nu))$ for some $\nu\in (0,2)$. There are some constants 
$\lambda_1,\lambda_2\in(0,\infty)$, depending only on $\nu,\kappa_1,\kappa_2$, such that
for all integer $n\geq 2$, all $v,\vs \in \rd$,
\begin{align*}
D_n(v,\vs):=&\int_0^\pi \int_0^{2\pi} [|v'|^{2n}+|v'_*|^{2n}-|v|^{2n}-|\vs|^{2n} ] \dd \varphi \beta(\theta)\dd\theta\\ 
\leq& -\lambda_1 n^{\nu/2} (|v|^{2n}+|\vs|^{2n}) \\
&+ \lambda_2
\sum_{a=1}^{n-1} \binom{n}{a} \Big(\frac{n^{\nu/2}}{(n-a)^{\nu/2+1}}+\frac 1a\Big)
 [|v|^{2a}|\vs|^{2(n-a)}+|v|^{2(n-a)}|\vs|^{2a}].
\end{align*}
\end{lem}

In the case with cutoff, see e.g. \cite{acgm}, one gets (roughly) something like
$$D_n(v,\vs) \leq -  (|v|^{2n}+|\vs|^{2n}) + \e_n
\sum_{a=1}^{n-1} \binom{n}{a} [|v|^{2a}|\vs|^{2(n-a)}+|v|^{2(n-a)}|\vs|^{2a}],
$$
with $\e_n\to 0$ as $n\to \infty$.
Here the negative term is reinforced by the factor $n^{\nu/2}$,
and this is the main fact we will have to exploit. We will also have to play tightly with
the positive term, showing that despite the fact it is not clearly multiplied by a small factor,
it can be absorbed, in some sense, by the negative term.

\vip

We start with an explicit computation of the $\varphi$-average.

\begin{lem}\label{groc}
For any integer $n\geq 2$, any $v,\vs \in \rd$, any $\theta \in (0,\pi]$, 
we have
\begin{align*}
\Theta_n(v,\vs,\theta):=\frac1{2\pi}\int_0^{2\pi} |v'|^{2n}\dd \varphi
= \Big(\frac{1+\cos\theta}2\Big)^n |v|^{2n} 
+\Big(\frac{1-\cos\theta}2\Big)^n |\vs|^{2n}+\Lambda_n(v,\vs,\theta),
\end{align*}
where, setting $\cA_n=\{(i,j,k)\in \nn^2 : i+j+k=n, \; i\leq n-1, \; j\leq n-1,\;k \in 2\nn\}$,
$$\Lambda_n(v,\vs,\theta)\!= \!\!\!\!\!\!
\sum_{(i,j,k)\in \cA_n} \!\!\! \frac{n!}{i!j! [(k/2)!]^2}\Big(\frac{1+\cos\theta}{2}\Big)^i
\Big(\frac{1-\cos\theta}{2}\Big)^j\Big(\frac{\sin \theta}{2}\Big)^k |v|^{2i}|\vs|^{2j}
\Big(|v|^2|\vs|^2-(v\cdot\vs)^2\Big)^{k/2}.
$$
\end{lem}

\begin{proof} We fix $n\geq 2$ and divide the proof into 3 steps.

\vip

{\it Step 1.} Recalling from \eqref{rep} that $v'=v-\frac{1-\cos\theta}2(v-\vs)+\frac{\sin\theta}2 \Gamma(v-\vs,\varphi)$,
that $| \Gamma(v-\vs,\varphi)|=|v-\vs|$
and that $(v-\vs)\cdot \Gamma(v-\vs,\varphi)=0$, we find
\begin{align*}
|v'|^2=&|v|^2\!+\!\Big(\frac{1\!-\!\cos\theta}2\Big)^2|v-\vs|^2\!+\!\Big(\frac{\sin\theta}2\Big)^2|v-\vs|^2
\!-\!(1-\cos\theta)v\cdot(v-\vs)\!+\!(\sin\theta)v\cdot\Gamma(v-\vs,\varphi)\\
=&\frac{1+\cos\theta}2 |v|^2 + \frac{1-\cos\theta}2|\vs|^2 + (\sin\theta) v\cdot\Gamma(v-\vs,\varphi).
\end{align*}
Hence, by Newton's trinomial expansion, setting $\cB_n=\{(i,j,k)\in\nn^3 : i+j+k=n\}$,
\begin{align*}
|v'|^{2n}=& \sum_{(i,j,k)\in \cB_n} \frac{n!}{i!j!k!}
\Big(\frac{1+\cos\theta}{2}\Big)^i\Big(\frac{1-\cos\theta}{2}\Big)^j(\sin \theta)^k 
|v|^{2i}|\vs|^{2j} (v\cdot\Gamma(v-\vs,\varphi))^k.
\end{align*}

{\it Step 2.} We now prove that for $k\in \nn$,
$$
\frac1{2\pi}\int_0^{2\pi} (v\cdot\Gamma(v-\vs,\varphi))^k \dd \varphi
= \indiq_{\{k \in 2\nn\}} \frac{k!}{2^k [(k/2)!]^2}(|v|^2|\vs|^2-(v\cdot\vs)^2)^{k/2}.
$$
We have $v\cdot\Gamma(v-\vs,\varphi)=a \cos\varphi + b \sin\varphi$, where $a=v\cdot I(v-\vs)$
and $b=v\cdot J(v-\vs)$, whence
$$
v\cdot\Gamma(v-\vs,\varphi)=\sqrt{a^2+b^2} \sin(\varphi+\varphi_0),
$$
for $\varphi_0$ such that $\frac{a}{\sqrt{a^2+b^2}}=\sin\varphi_0$
and $\frac{b}{\sqrt{a^2+b^2}}=\cos\varphi_0$.
We thus recognize a
Wallis integral:
$$
\frac1{2\pi}\int_0^{2\pi} (v\cdot\Gamma(v-\vs,\varphi))^k \dd \varphi
= \frac2\pi(a^2+b^2)^{k/2} \indiq_{\{k \in 2\nn\}} \int_0^{\pi/2} \sin^k \varphi\dd\varphi= 
\indiq_{\{k \in 2\nn\}} \frac{k!}{2^k [(k/2)!]^2}(a^2+b^2)^{k/2}.
$$
To complete the step, we recall that $(\frac{v-\vs}{|v-\vs|},\frac{I(v-\vs)}{|v-\vs|},
\frac{J(v-\vs)}{|v-\vs|})$ is an orthonormal basis, whence
$$
|v|^2= \frac{[v\cdot(v-\vs)]^2+[v\cdot I(v-\vs)]^2+[v\cdot J(v-\vs)]^2  }{|v-\vs|^2}=
\frac{[v\cdot(v-\vs)]^2+a^2+b^2}{|v-\vs|^2}
$$
and thus $a^2+b^2=|v|^2|v-\vs|^2-[v\cdot(v-\vs)]^2=|v|^2|\vs|^2-(v\cdot\vs)^2$.

\vip

{\it Step 3.} Gathering Steps 1 and 2 and setting $\cC_n=\{(i,j,k)\in\nn^3 : i+j+k=n,\; k\in 2\nn\}$,
\begin{align*}
\Theta_n(v,\vs,\theta)\!= \!\!\!\!\!\!
\sum_{(i,j,k)\in \cC_n} \!\!\! \frac{n!}{i!j! [(k/2)!]^2}\Big(\frac{1\!+\!\cos\theta}{2}\Big)^i
\Big(\frac{1\!-\!\cos\theta}{2}\Big)^j\Big(\frac{\sin \theta}{2}\Big)^k |v|^{2i}|\vs|^{2j}
\Big(|v|^2|\vs|^2\!-(v\cdot\vs)^2\Big)^{k/2}.
\end{align*}
It then suffices to isolate the two extreme terms $(i,j,k)=(n,0,0)$ and $(i,j,k)=(0,n,0)$.
\end{proof}

We next estimate, sharply, some integrals in $\theta$.

\begin{lem}\label{inte}
Assume $(H_2(\nu))$. There are $\zeta_1,\zeta_2 \in (0,\infty)$,
depending only on $\nu,\kappa_1,\kappa_2$, such that:
\vip
(i) for all integer $n\geq 2$,
$$
a_n:=\int_0^\pi \Big(1-\Big[\frac{1+\cos\theta}2\Big]^n-\Big[\frac{1-\cos\theta}2\Big]^n\Big) 
\beta(\theta) \dd \theta \geq \zeta_1 n^{\nu/2};
$$

(ii) for all integers $a,n$ such that $1\leq a \leq n-1$,
setting
$$
J_{n,a}:=\int_0^\pi \Big(\frac{1+\cos\theta}2\Big)^a\Big(\frac{1-\cos\theta}2\Big)^{n-a} \beta(\theta) \dd \theta
,\quad \hbox{we have}\quad \binom{n}{a} J_{n,a}
\leq \zeta_2 \Big[\frac{n^{\nu/2}}{(n-a)^{\nu/2+1}} + \frac{1}{a} \Big].
$$
\end{lem}

\begin{proof}
We start with (i): the integrand in $a_n$ is nonnegative, because
for $x=(1+\cos\theta)/2\in [0,1]$, we have $x^n+(1-x)^n\leq 1$. Hence recalling $(H_2(\nu))$,
$$
a_n \geq \kappa_1 \int_0^{n^{-1/2}}\Big(1-\Big[\frac{1+\cos\theta}2\Big]^n-\Big[\frac{1-\cos\theta}2\Big]^n\Big) 
\theta^{-\nu-1} \dd \theta.
$$
For all $\theta\in [0,n^{-1/2}]$, we have $(1+\cos \theta)/2\leq 1-\theta^2/5$,
whence $[(1+\cos \theta)/2]^n \leq (1-\theta^2/5)^n\leq \exp(-n \theta^2/5)$ and thus
$1-[(1+\cos \theta)/2]^n \geq n \theta^2/10$.
Next, still for $\theta\in [0,n^{-1/2}]$, we have  $(1-\cos \theta)/2\leq \theta^2/4$,
whence, since $n\geq 2$,
$$
1-\Big[\frac{1+\cos\theta}2\Big]^n-\Big[\frac{1-\cos\theta}2\Big]^n \geq 
\frac {n \theta^2}{10} -\frac{\theta^{2n}}{4^n} \geq \Big(\frac{n}{10}-\frac 1{16} \Big)\theta^2
\geq \Big(\frac{n}{10}-\frac n{32} \Big) \theta^2 \geq \frac{n}{20} \theta^2.
$$
Consequently,
$$
a_n \geq \kappa_1 \frac{n}{20} \int_0^{n^{-1/2}} \theta^{1-\nu} \dd \theta
= \frac{\kappa_1}{20(2-\nu)} n^{\nu/2}.
$$

\vip

For (ii), we first use $(H_2(\nu))$ to write $J_{n,a}\leq \kappa_2 K_{\nu,n,a}+\kappa_2 L_{\nu,n,a}$, where
\begin{align*}
K_{\nu,n,a}=& \int_0^{\pi/2} \Big(\frac{1+\cos\theta}2\Big)^a\Big(\frac{1-\cos\theta}2\Big)^{n-a} \theta^{-\nu-1}
\dd \theta,\\
L_{\nu,n,a} =& \Big(\frac 2\pi\Big)^{\nu+1}
\int_{\pi/2}^\pi \Big(\frac{1+\cos\theta}2\Big)^a\Big(\frac{1-\cos\theta}2\Big)^{n-a}\dd \theta.
\end{align*}
Using the substitution $\theta \to \pi-\theta$, we see that
\begin{align*}
L_{\nu,n,a}=&\Big(\frac 2\pi\Big)^{\nu+1}
\int_{0}^{\pi/2} \Big(\frac{1+\cos\theta}2\Big)^{n-a}\Big(\frac{1-\cos\theta}2\Big)^{a}\dd \theta\\
\leq& \Big(\frac 2\pi\Big)^{\nu}
\int_{0}^{\pi/2} \Big(\frac{1+\cos\theta}2\Big)^{n-a}\Big(\frac{1-\cos\theta}2\Big)^{a}\theta^{-1}\dd \theta
=\Big(\frac 2\pi\Big)^{\nu} K_{0,n,n-a}.
\end{align*}
We will prove that for $\nu\in[0,2)$, there is a constant $A_\nu\in (0,\infty)$ such that 
for all $1\leq a \leq n-1$,  
\begin{align}\label{gg}
\binom{n}{a} K_{\nu,n,a} \leq A_\nu \frac{n^{\nu/2}}{(n-a)^{1+\nu/2}}.
\end{align}
We will deduce that $\binom{n}{a} L_{\nu,n,a}\leq [2/\pi]^{\nu} \binom{n}{a} K_{0,n,n-a}\leq [2/\pi]^{\nu}A_0 a^{-1}$ 
and this will end the proof.
\vip
For $\theta\in (0,\pi/2]$, we have $\theta \leq 2\sin \theta$ and $\theta^{-1}\leq [(1-\cos\theta)/2]^{-1/2}$,
so that 
$$
\theta^{-\nu-1} \leq 2 \theta^{-\nu-2} \sin\theta \leq 2 \Big( \frac{1-\cos\theta}2\Big)^{-\nu/2-1} \sin\theta 
$$
and thus
$$
K_{\nu,n,a}\leq 2 \int_0^{\pi/2} \Big(\frac{1+\cos\theta}2\Big)^a\Big(\frac{1-\cos\theta}2\Big)^{n-a-\nu/2-1}
\sin\theta \dd \theta=4 \int_{1/2}^1 x^a (1-x)^{n-a-\nu/2-1}\dd x,
$$
using the change of variables $x=(1+\cos\theta)/2$. Hence
$$
 K_{\nu,n,a}\leq 4 
\int_0^{1} x^a (1-x)^{n-a-\nu/2-1}\dd x = 4 \frac{\Gamma(a+1)\Gamma(n-a-\nu/2)}{\Gamma(n+1-\nu/2)},
$$
where $\Gamma$ is Euler's Gamma function. Using that $\Gamma(k+1)=k!$, that
$(n-a-\nu/2)\Gamma(n-a-\nu/2)=\Gamma(n-a+1-\nu/2)$
and setting $u_{\nu,k}= \Gamma(k+1)/\Gamma(k+1-\nu/2)$, we realize that
$$
\binom{n}{a} K_{\nu,n,a} \leq 4 \frac{n!\Gamma(n-a-\nu/2)}{\Gamma(n+1-\nu/2)(n-a)!}
= 4 \frac{ u_{\nu,n}} {(n-a-\nu/2)u_{\nu,n-a}}.
$$
But using Stirling's formula $\Gamma(x+1)\sim \sqrt{2\pi x} (x/e)^x$ as $x\to \infty$, one can verify
that $u_{\nu,k}\sim k^{\nu/2}$ as $k\to \infty$ so that there is a constant $A_\nu \in (1,\infty)$ such that
for all $k\geq 1$,
$$
A_\nu^{-1} k^{\nu/2}\leq  u_{\nu,k} \leq A_\nu k^{\nu/2}.
$$
We conclude that for all $1\leq a \leq n-1$,
$$
\binom{n}{a} K_{\nu,n,a}\leq 4A_\nu^2 \frac{n^{\nu/2}}{(n-a)^{\nu/2}(n-a-\nu/2)}
\leq    \frac{4A_\nu^2}{1-\nu/2} \times\frac{n^{\nu/2}}{(n-a)^{\nu/2+1}}
$$
because $n-a\geq 1$ implies that $n-a-\nu/2 \geq (n-a)(1-\nu/2)$.
We have checked \eqref{gg} and the proof is complete.
\end{proof}

We can now handle the

\begin{proof}[Proof of Lemma \ref{pov}]
We fix $n\geq 2$. Using \eqref{rep}, we realize that, with the notation of Lemma \ref{groc},
\begin{align*}
\frac1{2\pi}\int_0^{2\pi} [|v'|^{2n}+|v'_*|^{2n}-|v|^{2n}-|\vs|^{2n} ] \dd \varphi
= & \Theta_n(v,\vs,\theta)+ \Theta_n(\vs,v,-\theta)-|v|^{2n}-|\vs|^{2n}.
\end{align*}
We deduce from Lemma \ref{groc} that
\begin{align*}
&\frac1{2\pi}\int_0^{2\pi} [|v'|^{2n}+|v'_*|^{2n}-|v|^{2n}-|\vs|^{2n} ] \dd \varphi\\
= & -\Big(1-\Big[\frac{1+\cos\theta}2\Big]^n-\Big[\frac{1-\cos\theta}2\Big]^n\Big)\times
 (|v|^{2n}+|\vs|^{2n})
  + \Lambda_n(v,\vs,\theta)+\Lambda_n(\vs,v,-\theta).
\end{align*}
Hence $D_n(v,\vs)=
-D_{n,1}(v,\vs)+D_{n,2}(v,\vs)$,
where
\begin{align*}
D_{n,1}(v,\vs)=& 2\pi\int_0^\pi \Big(1-\Big[\frac{1+\cos\theta}2\Big]^n-\Big[\frac{1-\cos\theta}2\Big]^n\Big) \beta(\theta) \dd \theta
\times (|v|^{2n}+|\vs|^{2n}),\\
D_{n,2}(v,\vs)=&2\pi \int_0^\pi [\Lambda_n(v,\vs,\theta)+\Lambda_n(\vs,v,-\theta)] \beta(\theta)\dd\theta.
\end{align*}
We now divide the proof into 5 steps.

\vip

{\it Step 1.}
By Lemma \ref{inte}-(i), we have $D_{n,1}(v,\vs)\geq 2\pi\zeta_1 n^{\nu/2}(|v|^{2n}+|\vs|^{2n})$.

\vip

{\it Step 2.}
We next roughly bound
$|v|^2|\vs|^2-(v\cdot\vs)^2$ by $|v|^2|\vs|^2$ in the expression of $\Lambda_n$ to find
$$
D_{n,2}(v,\vs)\leq 2\pi \sum_{(i,j,k)\in \cA_n} \!\!\! \frac{n!}{i!j! [(k/2)!]^2}I_{i,j,k}
[|v|^{2i+k}|\vs|^{2j+k}+|v|^{2j+k}|\vs|^{2i+k}],
$$
where 
$$
I_{i,j,k}=\int_0^{\pi} \Big(\frac{1+\cos\theta}{2}\Big)^i
\Big(\frac{1-\cos\theta}{2}\Big)^j\Big(\frac{\sin \theta}{2}\Big)^k \beta(\theta)\dd \theta.
$$
This can be rewritten as
\begin{align*}
D_{n,2}(v,\vs)\leq 2\pi \sum_{a=0}^n K_{n,a} [|v|^{2a}|\vs|^{2(n-a)}+|v|^{2(n-a)}|\vs|^{2a}]
\end{align*}
where, setting $\cA_{n,a}=\{(i,j,k)\in\cA_n : i+k/2=a$ (whence $j+k/2=n-a)\}$,
$$
K_{n,a}= \sum_{(i,j,k)\in \cA_{n,a}} \frac{n!}{i!j! [(k/2)!]^2} I_{i,j,k}.
$$

{\it Step 3.} We have $K_{n,n}=K_{n,0}=0$. Indeed, we e.g. have 
$\cA_{n,n}=\emptyset$ because for $(i,j,k)\in \cA_{n,n}$, we have $i+j+k=n$ and 
$j+k/2=0$, whence 
$i=n$, which is forbidden since $(i,j,k)\in\cA_n$.

\vip

{\it Step 4.} We now fix $a\in \ig 1,n-1 \id$. Then  $(i,j,k) \in \cA_{n,a}$ if and only if 
there is $\ell \in \nn$ such that 
$k=2\ell \in  \ig 0,n \id$, $i=a-\ell \in \ig 0,n-1\id$ and $j=(n-a)-\ell \in \ig 0,n-1 \id$,
so that
$$
K_{n,a}=\sum_{\ell=0}^{a\land(n-a)} \frac{n!}{(a-\ell)!(n-a-\ell)!(\ell!)^2} I_{a-\ell,n-a-\ell,2\ell}.
$$
But for any $\ell \in \ig 0,a\land(n-a)\id$,
it holds that 
$$
I_{a-\ell,n-a-\ell,2\ell}=\int_0^{\pi} \Big(\frac{1+\cos\theta}{2}\Big)^{a-\ell}
\Big(\frac{1-\cos\theta}{2}\Big)^{n-a-\ell}\Big(\frac{\sin \theta}{2}\Big)^{2\ell} \beta(\theta)\dd \theta
=J_{n,a},
$$
where
$$
J_{n,a}
=\int_0^{\pi} \Big(\frac{1+\cos\theta}{2}\Big)^{a}\Big(\frac{1-\cos\theta}{2}\Big)^{n-a}\beta(\theta)\dd \theta,
$$
because $(1+\cos\theta)(1-\cos\theta)=\sin^2\theta$. Thus
$$
K_{n,a}=J_{n,a}\sum_{\ell=0}^{a\land(n-a)} \frac{n!}{(a-\ell)!(n-a-\ell)!(\ell!)^2}
=J_{n,a} \binom{n}{a} \sum_{\ell=0}^{a\land(n-a)} \binom{a}{\ell} \binom{n-a}{\ell} = J_{n,a}\Big[\binom{n}{a}\Big]^2.
$$
We finally used the well-known identity $\sum_{\ell=0}^{a\land b}\binom{a}{\ell} \binom{b}{\ell}=
\binom{a+b}{a}$,
which can be shown noting that $\binom{a+b}{a}$ is the coefficient in front of $X^a$ of $(1+X)^{a+b}$,
while $\sum_{\ell=0}^{a\land b}\binom{a}{\ell} \binom{b}{\ell}=\sum_{\ell=0}^{a\land b}\binom{a}{\ell} \binom{b}{b-\ell}$ is the coefficient 
in front of $X^a$ of $(1+X)^{a}(1+X)^{b}$.

\vip

{\it Step 5.} Gathering Steps 2-3-4, we have checked that
\begin{align*}
D_{n,2}(v,\vs)\leq&2\pi \sum_{a=1}^{n-1} J_{n,a} \Big[ \binom{n}{a} \Big]^2 [|v|^{2a}|\vs|^{2(n-a)}+|v|^{2(n-a)}|\vs|^{2a}].
\end{align*}
The conclusion then follows from Lemma \ref{inte}-(ii), from which
$\binom{n}{a} J_{n,a}\leq \zeta_2 [\frac{n^{\nu/2}}{(n-a)^{\nu/2+1}} + \frac{1}{a}]$.
\end{proof}

\section{Even integer moments}\label{smom}

Using the previous Povzner inequality, we can derive the even integer moments.

\begin{lem}\label{2groc}
Assume  $(H_1(\gamma))$ and $(H_2(\nu))$ for some $\gamma \in (0,1]$ and some
$\nu \in (0,2)$. For any integer $n\geq 2$, any $t>0$, 
$$
m_{2n}'(t) \leq - c_1 n^{\nu/2} m_{2n+\gamma}(t)+ c_2 S_n(t) +  2c_1 n^{\nu/2}2^{2n/\gamma},
$$
where $c_1=\lambda_1$ and $c_2=16\lambda_2$ (see Lemma \ref{pov}) depend only on $\nu,\kappa_1,\kappa_2$
and where
$$
S_n(t) = \sum_{a=1}^{\lfloor n/2\rfloor} \binom{n}{a} \frac{n^{\nu/2}}{a^{\nu/2+1}} m_{2a}(t)m_{2(n-a)+\gamma}(t).
$$
\end{lem}

\begin{proof} We fix $n\geq 2$
and use the weak formulation \eqref{weak} with $\phi(v)=|v|^{2n}$, which is licit thanks to \eqref{momnul1}, 
to find,
for any $t>0$ and with the notation of Lemma \ref{pov},
$$
m_{2n}'(t)= \intrd\intrd D_n(v,\vs) |v-\vs|^\gamma f_t(\dd \vs) f_t(\dd v).
$$
Hence using Lemma \ref{pov}, $m_{2n}'(t)\leq -A_n(t)+B_n(t)$, where
\begin{align*}
A_n(t)=&\lambda_1 n^{\nu/2} 
\intrd \intrd (|v|^{2n}+|\vs|^{2n}) |v-\vs|^\gamma f_t(\dd \vs) f_t(\dd v),\\
B_n(t)=&\lambda_2 \sum_{a=1}^{n-1} \binom{n}{a}  \Big(\frac{n^{\nu/2}}{(n-a)^{\nu/2+1}}+\frac 1a\Big)
\intrd \intrd [|v|^{2a}|\vs|^{2(n-a)}+|v|^{2(n-a)}|\vs|^{2a}] \\
&\hskip9cm |v-\vs|^\gamma f_t(\dd \vs) f_t(\dd v).
\end{align*}
We now divide the proof into 2 steps.

\vip

{\it Step 1.} Using a symmetry argument and that $|v-\vs|^\gamma\geq |v|^\gamma-|\vs|^\gamma$, we see that
$$
A_n(t)=2\lambda_1 n^{\nu/2} 
\intrd \intrd |v|^{2n} |v-\vs|^\gamma f_t(\dd \vs) f_t(\dd v)\geq 2\lambda_1 n^{\nu/2} (m_{2n+\gamma}(t)-
m_{2n}(t)m_{\gamma}(t)).
$$
By H\"older's inequality, $m_\gamma(t)\leq [m_2(t)]^{\gamma/2}=1$. Using moreover that
$x^{2n} \leq \frac 12 x^{2n+\gamma} + 2^{2n/\gamma}$ for all $x\geq 0$ (separate the cases
$x \leq 2^{1/\gamma}$ and $x\geq 2^{1/\gamma}$),
$$
m_\gamma(t)m_{2n}(t)\leq m_{2n}(t)=\intrd |v|^{2n}f_t(\dd v) 
\leq \intrd \Big(\frac12|v|^{2n+\gamma}+2^{2n/\gamma}\Big)f_t(\dd v)=\frac12m_{2n+\gamma}(t) +  2^{2n/\gamma}.
$$
All in all,
$$
A_n(t)\geq  \lambda_1 n^{\nu/2} m_{2n+\gamma}(t) - 2\lambda_1  n^{\nu/2}2^{2n/\gamma}.
$$

\vip

{\it Step 2.} Using a symmetry argument and that $|v-\vs|^\gamma\leq |v|^\gamma+|\vs|^\gamma$, we find
\begin{align*}
B_n(t) \leq& 2\lambda_2\sum_{a=1}^{n-1} \binom{n}{a}  \Big(\frac{n^{\nu/2}}{(n-a)^{\nu/2+1}}+\frac 1a\Big)
[m_{2a+\gamma}(t)m_{2(n-a)}(t)+m_{2a}(t)m_{2(n-a)+\gamma}(t)]\\
=&2\lambda_2\sum_{a=1}^{n-1} \binom{n}{a}  \Big(\frac{n^{\nu/2}}{a^{\nu/2+1}}+\frac 1a\Big)
[m_{2a+\gamma}(t)m_{2(n-a)}(t)+m_{2a}(t)m_{2(n-a)+\gamma}(t)]
\end{align*}
by symmetry again. Since now $1/a \leq n^{\nu/2}/a^{\nu/2+1}$
and since $a\to n^{\nu/2}/a^{\nu/2+1}$ is decreasing,
$$
B_n(t) \leq 8\lambda_2 \sum_{a=1}^{\lfloor n/2\rfloor}\binom{n}{a} \frac{n^{\nu/2}}{a^{\nu/2+1}}
[m_{2a+\gamma}(t)m_{2(n-a)}(t)+m_{2a}(t)m_{2(n-a)+\gamma}(t)].
$$
But for $a\in \ig 1, \lfloor n/2\rfloor \id$, we have $a\leq n-a$ whence
$m_{2a+\gamma}(t)m_{2(n-a)}(t) \leq m_{2a}(t)m_{2(n-a)+\gamma}(t)$ by Lemma \ref{troc} below, so that
$$
B_n(t) \leq 16\lambda_2 \sum_{a=1}^{\lfloor n/2\rfloor}\binom{n}{a} \frac{n^{\nu/2}}{a^{\nu/2+1}}
m_{2a}(t)m_{2(n-a)+\gamma}(t)= 16\lambda_2  S_n(t).
$$
Recalling that $m_{2n}'(t)=-A_n(t)+B_n(t)$ and Step 1, we have reached our goal.
\end{proof}

We now prove a classical lemma that we used a few lines above.

\begin{lem}\label{troc}
For any $b\geq a \geq 0$, any $\alpha>0$, any $f\in\cP(\rd)$,
$$ 
m_{a+\alpha}(f)m_b(f)\leq
m_{b+\alpha}(f)m_a(f).
$$
\end{lem}

\begin{proof} We fix $\alpha>0$ and $f\in\cP(\rd)$ and have to prove that the function 
$u(a)=m_{a+\alpha}(f)/m_a(f)$ is nondecreasing on $\rr_+$. 
Observing that $\frac {\dd}{\dd a}m_a(f)=\intrd (\log |v|) |v|^a f(\dd v)$, we find
$$
u'(a)=\frac 1{m_a(f)}\intrd (\log |v|) |v|^{a+\alpha} f(\dd v)- \frac 1{m_a^2(f)}
\Big(\intrd (\log |v|) |v|^{a} f(\dd v)\Big)\Big(\intrd |v|^{a+\alpha} f(\dd v)\Big).
$$
Setting $g_a(\dd v)= |v|^{a} f(\dd v)/m_a(f)$, which is a probability measure,
$$
u'(a)=\frac 1\alpha \intrd (\log |v|^\alpha) |v|^{\alpha} g_a(\dd v)- \frac 1\alpha
\Big(\intrd (\log |v|^\alpha) g_a(\dd v)\Big)\Big(\intrd |v|^{\alpha} g_a(\dd v)\Big).
$$
Hence $u'(a)\geq 0$ by the Jensen inequality.
\end{proof}

We will also use the following remark.

\begin{rk}\label{hm2}
For $f \in \cP(\rd)$ such that $m_2(f)=1$, for all $r\geq s \geq 2$, by H\"older's inequality,
\begin{align*}
m_{s}(f) =\intrd |v|^{s-2} |v|^2 f(\dd v) \leq \Big(\intrd |v|^{r-2} |v|^2 f(\dd v)\Big)^{(s-2)/(r-2)}
=[m_r(f)]^{(s-2)/(r-2)}.
\end{align*}
\end{rk}

We finally quickly prove some moment estimates that are more or less well-known.

\begin{lem}\label{mm}
Assume $(H_1(\gamma))$ and $(H_2(\nu))$ for some $\gamma \in (0,1]$ and some
$\nu \in (0,2)$.
\vip
(i) For all $r>0$, there is $K_r\in (0,\infty)$, depending only on $\gamma,\nu,\kappa_1,\kappa_2,r$,
such that for all $t>0$, $$m_r(t)\leq K_r(1+t^{-(r-2)/\gamma}).$$

(ii) For all integer $n\geq 2$, all $A\geq 1$, there is $K_{2n,A}\in (0,\infty)$, depending only on 
$\gamma,\nu,\kappa_1,\kappa_2,n,A$,
{such that} 
$$ m_{2n}(0)\leq A \quad \Longrightarrow \quad \sup_{t\geq 0}m_{2n}(t)\leq K_{2n,A}.$$
\end{lem}

\begin{proof}
We first prove that for any fixed integer $n\geq 2$, there is a constant $A_n$, depending only on 
$\gamma,\nu,\kappa_1,\kappa_2,n$, allowed to vary from line to line, such that for all $t>0$,
\begin{equation}\label{cqff}
m_{2n}'(t) \leq -\frac{c_1}2 [m_{2n}(t)]^{1+\gamma/(2n-2)}+A_n.
\end{equation}
Using Lemma \ref{2groc}, a rough upperbound, and then Lemma \ref{troc}, we get
$$
m_{2n}'(t) \!\leq\! -c_1 m_{2n+\gamma}(t) + A_n \! \sum_{a=1}^{\lfloor n/2\rfloor} \!m_{2a}(t)m_{2(n-a)+\gamma}(t)
+A_n
\! \leq \! -c_1 m_{2n+\gamma}(t) +A_nm_{2(n-1)+\gamma}m_2(t) +A_n.
$$
Using now Remark \ref{hm2} and that $m_2(t)=1$, we find
\begin{align*}
m_{2n}'(t)\leq   - c_1[m_{2n}(t)]^{1+\gamma/(2n-2)} + A_n [m_{2n}(t)]^{1-(2-\gamma)/(2n-2)} +A_n,
\end{align*}
from which \eqref{cqff} follows. 
\vip
Point (ii) clearly follows from \eqref{2groc}. When $r=2n\geq 4$ is an even integer, point (i) also
follows from \eqref{2groc}. For a general $r>2$, we consider an integer $n_r\geq 2$ such that 
$2n_r\geq r$ and we conclude from Remark \ref{hm2} that 
$$m_r(t)\leq [m_{2n_r}(t)]^{(r-2)/(2n_r-2)}
\leq (K_{2n_r}(1+t^{-(2n_r-2)/\gamma}))^{(r-2)/(2n_r-2)},$$ 
whence the result. 
Finally, for $r\in (0,2]$, we obviously have $m_r(t)\leq [m_2(t)]^{r/2}\leq 1$.
\end{proof}

\section{Control of exponential moments by even integer moments}\label{sr}

We explain how to control exponential moments from even integer moments and {\it vice versa}.

\begin{lem}\label{biznorm} 
(i) For $f \in \cP(\rd)$, for $\sigma_0\in(0,\infty)$, $\alpha\in [1,\infty)$ and
$K \in [1,\infty)$,
$$\sup_{n\geq 0} \frac{\sigma_0 ^nm_{2n}(f)}{ (n!)^\alpha} \leq K
\quad \Longrightarrow \quad 
\intrd \exp(\sigma_0^{1/\alpha}|v|^{2/\alpha}/2) f(\dd v) \leq 2K^{1/\alpha}.
$$
(ii) For $\rho\in (0,2]$, $\sigma_0\in(0,1]$ and $K\in(1,\infty)$, 
there is $\sigma_1\in(0,\infty)$, depending only on $\rho,\sigma_0,K$, such that for all $f \in \cP(\rd)$,
$$
\intrd \exp(\sigma_0|v|^{\rho}) f(\dd v) \leq K \quad \Longrightarrow \quad
\sup_{n\geq 0} \frac{\sigma_1^nm_{2n}(f)}{ (n!)^{2/\rho}} \leq 1.
$$
\end{lem}

\begin{proof}
We start with (i). By H\"older's inequality, since $\alpha\geq 1$,
$$
m_{2n/\alpha}(f) \leq [m_{2n}(f)]^{1/\alpha}\leq K^{1/\alpha} \sigma_0^{-n/\alpha} n!,
$$
whence
$$
\intrd \exp[\sigma_0^{1/\alpha}|v|^{2/\alpha}/2] f(\dd v)=\sum_{n\geq 0} \frac{\sigma_0^{n/\alpha} m_{2n/\alpha}(f) }{2^nn!}
\leq K^{1/\alpha}  \sum_{n\geq 0} 2^{-n}=2K^{1/\alpha}.
$$

We now check (ii). We know that 
$$
\sup_{n\geq 0} \frac{\sigma_0^n m_{\rho n}(f)}{n!}\leq
\sum_{n\geq 0} \frac{\sigma_0^n m_{\rho n}(f)}{n!}=\intrd \exp(\sigma_0|v|^{\rho}) f(\dd v)\leq K.
$$
For $n\geq 1$, we set $k_n=\lceil 2n/\rho \rceil \in [2n/\rho,2n/\rho+1)$ and write,
using that $\rho k_n \geq 2n$,
$$
m_{2n}(f) \leq 1+ m_{\rho k_n}(f) \leq 1+\frac{K k_n!}{\sigma_0^{k_n}}, \quad 
\hbox{so that}\quad  \frac{\sigma_1^nm_{2n}(f)}{(n!)^{2/\rho}} \leq 
\frac{\sigma_1^n}{(n!)^{2/\rho}}+
\frac{\sigma_1^n K k_n!}{\sigma_0^{k_n} 
(n!)^{2/\rho}}.
$$
By Stirling's formula $n!\sim \sqrt{2\pi n} (n/e)^n$ and since $k_n=\lceil 2n/\rho \rceil$, we find that, 
for some constant $A\in (0,\infty)$, depending only on $\rho$ and allowed to vary, for all $n\geq 1$,
$$
\frac{k_n!}{(n!)^{2/\rho}} \!\leq\! A \frac{n^{1/2} [(2n/\rho+1)/e]^{2n/\rho+1}}{n^{1/\rho}[n/e]^{2n/\rho}}
\!\leq\! A n^{3/2-1/\rho} \frac{[(3n)/(e\rho)]^{2n/\rho}}{[n/e]^{2n/\rho}}\!\leq\!A n^{3/2}\Big(\frac 3 \rho \Big)^{2n/\rho}
\!\!\! \leq A \Big(\frac 4 \rho \Big)^{2n/\rho}.
$$
Observing next that $\sigma_0^{k_n} \geq \sigma_0^{2n/\rho+1}$ (since $\sigma_0\in (0,1]$), we end with
$$
\hbox{for all $n\geq 1$,}\quad 
\frac{\sigma_1^n m_{2n}(f)}{(n!)^{2/\rho}} \leq \frac{\sigma_1^n}{(n!)^{2/\rho}}+\frac {KA} {\sigma_0} \frac{\sigma_1^n}
{\sigma_0^{2n/\rho}} \Big(\frac 4 \rho \Big)^{2n/\rho}.
$$
This last quantity is bounded by $1$ if $\sigma_1>0$ is small enough (depending on
$\rho$, $K$, and $\sigma_0$).
\end{proof}

\section{Creation of exponential moments}\label{sc}

The following estimate will allow us to prove the creation result by Lemma \ref{biznorm}.

\begin{lem}\label{crea1}
Assume $(H_1(\gamma))$ and $(H_2(\nu))$ for some $\gamma \in (0,1]$ and some
$\nu \in (0,2)$ and set $\alpha =\max\{1,(2-\nu)/\gamma\}$.
There are $\sigma\in (0,1]$ and $T>0$, depending only on $\gamma,\nu,\kappa_1,\kappa_2$,
such that 
$$
\sup_{t\in [0,T]} \sum_{n=0}^\infty  \frac{(\sigma t)^{2n/\gamma}m_{2n}(t) }{ (n!)^\alpha} \leq 2.
$$
\end{lem}

\begin{proof} 
We recall that by Lemma \ref{2groc}, for any integer $n\geq 2$, any $t>0$,
\begin{equation}\label{rec}
m_{2n}'(t) \leq - c_1 n^{\nu/2} m_{2n+\gamma}(t)+ c_2 S_n(t) +  2c_1 n^{\nu/2}2^{2n/\gamma},
\end{equation}
where $S_n(t) = \sum_{a=1}^{\lfloor n/2\rfloor} \binom{n}{a} \frac{n^{\nu/2}}{a^{\nu/2+1}} m_{2a}(t)m_{2(n-a)+\gamma}(t)$.

\vip

{\it Step 1.} We introduce,
for $\sigma \in (0,1]$ to be chosen later, for $p\geq 2$
and $t\geq 0$,
$$
E_p(t)=\sum_{n=0}^p  \frac{(\sigma t)^{2n/\gamma}m_{2n}(t)}{(n!)^\alpha}.
$$
By Lemma \ref{mm}-(i) and since $m_0(t)=1$, it holds that for some constant $C_p\in (0,\infty)$,
$$
1\leq E_p(t) \leq 1 + C_p \sum_{n=1}^p t^{2n/\gamma}m_{2n}(t)\leq 1+C_p\sum_{n=1}^p t^{2n/\gamma}[1+t^{-(2n-2)/\gamma}]
=1+C_p\sum_{n=1}^p [t^{2n/\gamma}+t^{2/\gamma}],
$$
whence $\lim_{t\to 0} E_p(t)=1$.

\vip

{\it Step 2.} Since $m_0'(t)=m_{2}'(t)=0$, we deduce from \eqref{rec} that for all $p\geq 2$, all $t\in [0,1]$,
$$
E_p'(t) \leq -c_1 F_p(t) + c_2 G_p(t)+\frac{2\sigma}\gamma H_p(t) + C,
$$
where
\begin{gather*}
F_p(t)=\sum_{n=2}^p n^{\nu/2} \frac{(\sigma t)^{2n/\gamma}m_{2n+\gamma}(t) }{(n!)^\alpha},\quad
G_p(t)=\sum_{n=2}^p   \frac{(\sigma t)^{2n/\gamma}S_n(t)}{(n!)^\alpha},\quad
H_p(t)=\sum_{n=1}^p \frac{n (\sigma t)^{2n/\gamma -1} m_{2n}(t)}{(n!)^\alpha}
\end{gather*}
and where, since $\sigma t\leq 1$,
$$
C= 2 c_1 \sum_{n=2}^\infty  \frac{n^{\nu/2}2^{2n/\gamma}}{(n!)^\alpha}<\infty.
$$

{\it Step 3.} We first prove that for all $\e\in(0,\infty)$, there is  $A_\e\in(0,\infty)$,
depending only on $\gamma,\nu,\kappa_1,\kappa_2,\e$,
such that for any choice of $\sigma \in (0,1]$, for all $p\geq 2$, all $t\in[0,1]$,
$$
G_p(t) \leq \e E_p(t)F_p(t) + \sigma^{2/\gamma} A_\e (F_p(t)+1).
$$
We start from 
\begin{align*}
G_p(t)=&\sum_{n=2}^p \sum_{a=1}^{\lfloor n/2\rfloor} \frac{(\sigma t)^{2n/\gamma}}{(n!)^\alpha}\frac{n^{\nu/2}}{a^{\nu/2+1}}
\binom{n}{a} m_{2a}(t)m_{2(n-a)+\gamma}(t).
\end{align*}
Since $\alpha \geq 1$ and since $a \leq\lfloor n/2\rfloor$ implies $n \leq 2(n-a)$
and thus $n^{\nu/2} \leq 2^{\nu/2}(n-a)^{\nu/2}$,
\begin{align*}
G_p(t)
\leq& 2^{\nu/2}  \sum_{n=2}^p \sum_{a=1}^{\lfloor n/2\rfloor} 
\frac{(\sigma t)^{2n/\gamma}}{(n!)^\alpha}\frac{(n-a)^{\nu/2}}{a^{\nu/2+1}}
\Big[\binom{n}{a}\Big]^\alpha m_{2a}(t)m_{2(n-a)+\gamma}(t)\\
=& 2^{\nu/2}  \sum_{n=2}^p \sum_{a=1}^{\lfloor n/2\rfloor} (n-a)^{\nu/2}
\frac{(\sigma t)^{2(n-a)/\gamma}m_{2(n-a)+\gamma}(t)}{((n-a)!)^\alpha}  
\frac{(\sigma t)^{2a/\gamma}m_{2a}(t)}{(a!)^\alpha a^{\nu/2+1}}\\
=& 2^{\nu/2}  \sum_{a=1}^{\lfloor p/2\rfloor} \sum_{n=2a}^{p} (n-a)^{\nu/2}
\frac{(\sigma t)^{2(n-a)/\gamma}m_{2(n-a)+\gamma}(t)}{((n-a)!)^\alpha}  
\frac{(\sigma t)^{2a/\gamma}m_{2a}(t)}{(a!)^\alpha a^{\nu/2+1}}.
\end{align*}
Using the change of indices $n \to \ell=n-a$
\begin{align*}
G_p(t) \leq & 2^{\nu/2}  \sum_{a=1}^{\lfloor p/2\rfloor} \sum_{\ell=a}^{p-a} \ell^{\nu/2}
\frac{(\sigma t)^{2\ell/\gamma}m_{2\ell+\gamma}(t)}{(\ell!)^\alpha}  
\frac{(\sigma t)^{2a/\gamma}m_{2a}(t)}{(a!)^\alpha a^{\nu/2+1}}\\
\leq& 
 \Big(\sum_{\ell=1}^{p} \ell^{\nu/2}
\frac{(\sigma t)^{2\ell/\gamma}m_{2\ell+\gamma}(t)}{(\ell!)^\alpha}\Big) \times 2^{\nu/2}
\Big(\sum_{a=1}^p  \frac{(\sigma t)^{2a/\gamma}m_{2a}(t)}{(a!)^\alpha a^{\nu/2+1}} \Big)\\
=& \Big(F_p(t)+ (\sigma t)^{2/\gamma}m_{2+\gamma}(t)  \Big) I_p(t),
\end{align*}
where 
\begin{align*}
I_p(t)=  2^{\nu/2} \sum_{a=1}^p  \frac{(\sigma t)^{2a/\gamma}m_{2a}(t)}{(a!)^\alpha a^{\nu/2+1}} .
\end{align*}
But setting $N_\e= \lceil 2^{\nu/(\nu+2)} \e^{-2/(\nu+2)}\rceil$, it holds that
$a\geq N_\e$ implies $2^{\nu/2}a^{-\nu/2-1} \leq \e$, whence
$$
I_p(t)\leq \e E_p(t)+ J_{\e}(t),
\quad \hbox{where}\quad J_{\e}(t)=2^{\nu/2} \sum_{a=1}^{N_\e}
\frac{(\sigma t)^{2a/\gamma}m_{2a}(t)}{(a!)^\alpha}.
$$
By Lemma \ref{mm}-(i), we see that for some constants $A,A_\e\in(0,\infty)$, depending only on
$\gamma,\nu,\kappa_1,\kappa_2,\e$ and allowed to vary,
for any choice of $\sigma\in (0,1]$,
for all $t\in [0,1]$,
$$
(\sigma t)^{2/\gamma}m_{2+\gamma}(t)\leq A (\sigma t)^{2/\gamma}(1+t^{-1})\leq A \sigma^{2/\gamma} 
\;\; \hbox{and}\;\;
J_\e(t) \leq A_\e \sum_{a=1}^{N_\e} (\sigma t)^{2a/\gamma} [1+t^{-(2a-2)/\gamma}]
\leq A_\e \sigma^{2/\gamma}.
$$
All in all, we have proved that for any choice of $\sigma\in (0,1]$,
for all $t\in [0,1]$,
$$
G_p(t) \leq (F_p(t)+A\sigma^{2/\gamma} )(\e E_p(t) + A_\e \sigma^{2/\gamma} )
\leq \e E_p(t)F_p(t)+ \sigma^{2/\gamma}[\e A E_p(t)+ A_\e F_p(t)+AA_\e ].
$$
Since $E_p(t)\leq m_0(t)+ \sigma t m_{2}(t)+F_p(t)\leq 2+F_p(t)$, we conclude that, as desired,
$$
G_p(t) \leq \e E_p(t)F_p(t)+ \sigma^{2/\gamma}[(\e A+A_\e)F_p(t)+AA_\e +2\e A].
$$

\vip

{\it Step 4.} We now verify that there are some constants $\kappa, B\in(0,\infty)$,
depending only on $\gamma,\nu,\kappa_1,\kappa_2$, such that for any choice
of $\sigma \in (0,1]$, for all $p\geq 2$, all $t\in [0,1]$,
$$H_p(t) \leq \kappa F_p(t)+B.$$

We first observe that, for $\kappa>0$ to be chosen later,  
$$
\frac{n m_{2n}(t)}{\sigma t} \leq \kappa n^{\nu/2} [m_{2n}(t)]^{1+\gamma/(2n-2)} + \frac{n}{\sigma t}\Big(
\frac {n^{1-{\nu/2}}}{\kappa\sigma t}\Big)^{(2n-2)/\gamma}.
$$
Indeed, $nm_{2n}(t)/(\sigma t)$ is bounded by the second term if 
$m_{2n}(t) \leq [n^{1-{\nu/2}}/(\kappa\sigma t)]^{(2n-2)/\gamma}$
and by the first term else. Since $[m_{2n}(t)]^{1+\gamma/(2n-2)}\leq m_{2n+\gamma}(t)$ by Remark \ref{hm2}, we conclude
that
\begin{align*}{}
H_p(t) =&\sum_{n=1}^p \frac{n m_{2n}(t)}{\sigma t} \frac{(\sigma t)^{2n/\gamma}}{(n!)^\alpha}\\
\leq & \kappa \sum_{n=1}^p n^{\nu/2} \frac{(\sigma t)^{2n/\gamma}m_{2n+\gamma}(t) }{(n!)^\alpha}
+ \sum_{n=1}^p \frac{n}{\sigma t} \Big(
\frac {n^{1-{\nu/2}}}{\kappa \sigma t}\Big)^{(2n-2)/\gamma} \frac{(\sigma t)^{2n/\gamma}}{(n!)^\alpha}\\
=&\kappa F_p(t)+ \kappa (\sigma t)^{2/\gamma}m_{2+\gamma}(t)+(\sigma t)^{2/\gamma-1}\sum_{n=1}^p 
\frac{n^{(2-\nu)n / \gamma - (2-\nu)/ \gamma  +1}  }
{(n!)^\alpha \kappa^{(2n-2)/\gamma}}.
\end{align*}
But by Lemma \ref{mm}-(i), there is $A\in (0,\infty)$, depending only on $\gamma,\nu,\kappa_1,\kappa_2$,
such that for all $t \in [0,1]$,
$$
t^{2/\gamma}m_{2+\gamma}(t) \leq A t^{2/\gamma}[1+t^{-1}] \leq 2A. 
$$
Moreover, since $\alpha=\max\{1,(2-\nu)/\gamma\}\geq (2-\nu)/\gamma$, the series
$$
S:=\sum_{n=1}^\infty \frac{n^{(2-\nu)n / \gamma - (2-\nu)/ \gamma  +1}  }{(n!)^\alpha\kappa^{(2n-2)/\gamma}}
\leq \kappa^{2/\gamma} \sum_{n=1}^\infty \frac{n^{\alpha n +1}  }{(n!)^{\alpha}\kappa^{2n/\gamma}}
$$
is convergent if $\kappa=2e^{\alpha\gamma/2}$, by Stirling's formula $n!\sim \sqrt{2\pi n}(n/e)^n$.
Hence for any choice of $\sigma \in (0,1]$, for all $t\in [0,1]$,
$$
H_p(t) \leq \kappa F_p(t) + 2\kappa A + S.
$$

{\it Step 5.} By Steps 2-3-4, for any choice of $\e\in(0,\infty)$ and $\sigma \in (0,1]$,  
for all $p\geq 2$, all  $t\in (0,1]$,
\begin{align*}
E_p'(t) \leq & -c_1 F_p(t)+ c_2 G_p(t)+\frac{2\sigma}\gamma  H_p(t)+C\\
\leq & -c_1F_p(t)+c_2 [\e E_p(t)F_p(t) + \sigma^{2/\gamma} A_\e (F_p(t)+1)] 
+ \frac{2\sigma}\gamma [\kappa F_p(t)+B] + C.
\end{align*}
Choosing $\e=c_1/(4c_2)$ and $\sigma \in (0,1]$ small enough so that
$c_2 \sigma^{2/\gamma} A_\e + 2\kappa \sigma  /\gamma \leq c_1/2$, we find
$$
E_p'(t) \leq  -\frac{c_1}2 F_p(t) + \frac{c_1}4 E_p(t)F_p(t) + D,
$$
for some constant $D\in (0,\infty)$.
By Step 1, we know that $E_p(0)=\lim_{t\to 0}E_p(t)=1$, whence
$$
T_p = \sup \{ t>0 : E_p(t)\leq 2\} >0.
$$
For $p\geq 2$ and $t \in [0,T_p\land 1)$, we have $E_p'(t) \leq D$. Recalling that $E_p(0)=1$, we deduce that 
$T_p\geq \min\{1,1/D\}=:T.$ 
In other words, with our choice of $\sigma$,
all $t\in [0,T]$, $\sup_{p\geq 2}E_p(t) \leq 2$.
\end{proof}

We can finally give the

\begin{proof}[Proof of Theorem \ref{mr}-(i)]
We assume $(H_1(\gamma))$ and $(H_2(\nu))$ for some $\gamma \in (0,1]$ and some $\nu\in (0,2)$.
We fix $\rho=\min\{2\gamma/(2-\nu),2\}$ and observe that
$\alpha:=2/\rho = \max\{1,(2-\nu)/\gamma\}$. Thus by Lemma \ref{crea1}, 
there are $\sigma\in (0,1]$ and $T>0$, depending only on $\gamma,\nu,\kappa_1,\kappa_2$, such that
$$
\sup_{[0,T]} \sum_{n=0}^\infty \frac{(\sigma t)^{2n/\gamma}m_{2n}(t) }{ (n!)^\alpha} \leq 2.
$$
We deduce from Lemma \ref{biznorm}-(i) (with $\sigma_0=(\sigma t)^{2/\gamma}$) that
$$
\sup_{[0,T]} \intrd \exp[(\sigma t)^{\rho/\gamma} |v|^{\rho}/2] f_t(\dd v)
=
\sup_{[0,T]} \intrd \exp[(\sigma t)^{2/(\gamma\alpha)} |v|^{2/\alpha}/2] f_t(\dd v )
\leq 2^{1+1/\alpha}\leq 4.
$$
as desired.
\end{proof}

\section{Propagation of exponential moments}\label{sp}

We proceed as in the previous section.

\begin{lem}\label{propa1}
Assume $(H_1(\gamma))$ and $(H_2(\gamma))$ for some $\gamma \in (0,1]$ and some $\nu\in (0,2)$.
For any $\sigma_0 \in(0,\infty)$ and any $\alpha \geq 1$, there
is  $\sigma \in (0,\sigma_0]$, depending only on $\gamma,\nu,\kappa_1,\kappa_2,\sigma_0$,
such that
\begin{align*}
\sup_{n\geq 0} \frac{\sigma_0^nm_{2n}(0)}{(n!)^\alpha}\leq 1
\quad \Longrightarrow \quad \sup_{t\geq 0} \sum_{n=0}^\infty \frac{\sigma^n m_{2n}(t)}{ (n!)^\alpha}\leq 3.
\end{align*}
\end{lem}

\begin{proof}  We fix $\alpha \geq 1$ and $\sigma_0\in (0,\infty)$ 
and assume that $\sup_{n\geq 0} (n!)^{-\alpha} \sigma_0^nm_{2n}(0)\leq 1$.

\vip 

{\it Step 1.} We introduce, for $\sigma \in (0,1]$ to be chosen later, for $p\geq 2$
and $t\geq 0$,
$$
E_p(t)=\sum_{n=0}^p \frac{\sigma^n m_{2n}(t)}{(n!)^\alpha}.
$$
If $\sigma \in (0,\sigma_0/2]$, 
$E_p(0) \leq 2$ for all $p\geq 2$, because
$$
E_p(0)\leq 1+\sum_{n\geq 1} \frac{(\sigma_0/2)^n m_{2n}(0)}{(n!)^\alpha} 
\leq 1+\sum_{n\geq 1} 2^{-n}=2.
$$

{\it Step 2.} By Lemma \ref{2groc}, since $m_0'(t)=m_2'(t)=0$ and since $\sigma\in (0,1]$,
$$
E_p'(t) \leq -c_1 F_p(t) + c_2 G_p(t)+ C,
$$
where
\begin{gather*}
F_p(t)=\sum_{n=2}^p n^{\nu/2} \frac{\sigma^n m_{2n+\gamma}(t) }{(n!)^\alpha},\qquad
G_p(t)=\sum_{n=2}^p \frac{\sigma^nS_n(t)}{(n!)^\alpha}, \qquad
C= 2c_1 \sum_{n=2}^\infty  \frac{n^{\nu/2}2^{2n/\gamma}}{(n!)^\alpha}.
\end{gather*}

{\it Step 3.} We prove that for all $\e\in(0,\infty)$, there is a constant  $A_\e\in(0,\infty)$,
depending only on $\e,\gamma,\nu,\kappa_1,\kappa_2,\alpha,\sigma_0$,
such that for any choice of $\sigma \in (0,1]$, for all $p\geq 2$, all $t\geq 0$,
$$
G_p(t) \leq \e E_p(t)F_p(t) + \sigma A_\e (F_p(t)+1).
$$
Exactly as in the previous section, Step 3, we have
\begin{align*}
G_p(t) \leq& \Big(F_p(t)+ \sigma m_{2+\gamma}(t)  \Big) \Big(\e E_p(t)+ J_{\e}(t) \Big),
\end{align*}
where, setting $N_\e= \lceil 2^{\nu/(\nu+2)} \e^{-2/(\nu+2)}\rceil$,
\begin{align*}
J_{\e}(t)=2^{\nu/2} \sum_{a=1}^{N_\e}
\frac{\sigma^{a}m_{2a}(t)}{(a!)^\alpha},
\end{align*}
For all $n\geq 2$, $m_{2n}(0)\leq \sigma_0^{-n}(n!)^\alpha$, so that
by Lemma \ref{mm}-(ii), there are some constants 
$A,A_\e\in(0,\infty)$, 
depending only on $\e,\gamma,\nu,\kappa_1,\kappa_2,\alpha,\sigma_0$ and allowed to vary, such that
for all $t\geq 0$,
$$
m_{2+\gamma}(t)\leq [m_4(t)]^{(2+\gamma)/4}\leq A \quad \hbox{and}\quad
J_\e(t) \leq A_\e \sum_{a=1}^{N_\e} \sigma^{a} \leq A_\e \sigma.
$$
All in all, we have proved that for any choice of $\sigma\in (0,1]$,
for all $t\geq 0$,
$$
G_p(t) \leq (F_p(t)+A \sigma )(\e E_p(t) + A_\e \sigma )
\leq \e E_p(t)F_p(t)+ \sigma[\e A E_p(t)+ A_\e F_p(t)+AA_\e ].
$$
The conclusion follows since $E_p(t)\leq m_0(t)+ \sigma m_{2}(t)+F_p(t)\leq 2+F_p(t)$.

\vip

{\it Step 4.} We now prove that for all $p\geq 2$, all $t\geq 0$,
$$F_p(t)\geq \sigma^{-\gamma/2}(E_p(t)-e).$$
We write
\begin{align*}
F_p(t)=&\sum_{n=2}^p n^{\nu/2} \frac{\sigma^nm_{2n+\gamma}(t)}{(n!)^\alpha}\geq \frac{1}{\sigma^{\gamma/2}}
\sum_{n=2}^p \frac{\sigma^{n+\gamma/2}m_{2n+\gamma}(t)}{(n!)^\alpha}\geq\frac{1}{\sigma^{\gamma/2}}
\sum_{n=2}^p \frac{\sigma^{n}m_{2n}(t) -1}{(n!)^\alpha},
\end{align*}
because
$$
\sigma^{n+\gamma/2}m_{2n+\gamma}(t)=\intrd (\sigma |v|^2)^{n+\gamma/2}f_t(\dd v) \geq 
\intrd [(\sigma |v|^2)^{n}-1]f_t(\dd v) = \sigma^{n}m_{2n}(t) -1.
$$
Hence, since  $m_0(t)=m_2(t)=1$, since $\sigma \in (0,1]$ and since $\alpha\geq 1$,
$$
F_p(t)\geq\frac{1}{\sigma^{\gamma/2}}\Big[E_p(t)-m_0(t)-\sigma m_{2}(t)
-\sum_{n=2}^p \frac{1}{(n!)^\alpha}\Big]
\geq \frac{1}{\sigma^{\gamma/2}}\Big[E_p(t)-\sum_{n=0}^p \frac{1}{(n!)^\alpha}\Big]\geq
\frac{1}{\sigma^{\gamma/2}}[E_p(t)-e].
$$

{\it Step 5.} By Steps 2-3-4, for any choice of $\e\in(0,\infty)$ and $\sigma \in (0,1]$,  
for all $p\geq 2$, all  $t\geq 0$,
\begin{align*}
E_p'(t) \leq & -c_1 F_p(t) + c_2 G_p(t)+C \leq
-c_1 F_p(t) + c_2[\e E_p(t)F_p(t)+ \sigma A_\e(F_p(t)+1)]+C.
\end{align*}
Choosing $\e=c_1/(16c_2)$ and then $\sigma_1 \in (0,1]$ small enough so that
$c_2 \sigma_1 A_\e \leq c_1/2$, 
we conclude that for some constant $D\in(0,\infty)$, for any choice of
$\sigma \in (0,\sigma_1]$, any $t\geq 0$,
$$
E_p'(t) \leq  -\frac{c_1}2 F_p(t) + \frac{c_1}{16} E_p(t)F_p(t) +D.
$$
We now recall from Step 1 that for $\sigma \in (0,\sigma_0/2]$, we have 
$E_p(0)\leq 2$ and by continuity,
$$
T_p= \sup\{t\geq 0 : E_p(t)\leq 4\}>0.
$$
But for all $t \in [0,T_p)$, if $\sigma \in (0,(\sigma_0/2)\land\sigma_1]$,
$$
E_p'(t) \leq -\frac{c_1}{4}F_p(t) + D \leq -\frac{c_1}{4\sigma^{\gamma/2}}[E_p(t)-e]+D
$$
by Step 4. Since $E_p(0)\leq 2$, this implies that  if
$\sigma \in (0,(\sigma_0/2)\land\sigma_1]$, for all $p\geq 2$ and all $t\in [0,T_p)$,
$$
E_p(t) \leq e+ \frac{4\sigma^{\gamma/2}D}{c_1}.
$$
Choosing $\sigma \in (0,(\sigma_0/2)\land\sigma_1]$ small enough so that ${4\sigma^{\gamma/2}D}/{c_1}\leq 3-e$,
we conclude that for all $p\geq 2$, $E_p(t)\leq 3$ for all $t\in[0,T_p)$, whence $T_p=\infty$
by continuity. In other words, for all $p\geq 2$, all $t\geq 0$, $E_p(t)\leq 3$, which was our goal.
\end{proof}

We finally can give the

\begin{proof}[Proof of Theorem \ref{mr}-(ii)]
We assume $(H_1(\gamma))$ and $(H_2(\nu))$ for some $\gamma \in (0,1]$ and some $\nu\in (0,2)$.
We fix $A>1$, $\rho \in (0,2]$ and $\sigma_0>0$ and assume that 
$\intrd \exp(\sigma_0 |v|^\rho)f_0(\dd v)\leq A$. 

\vip

We set $\alpha=2/\rho\geq 1$. By Lemma \ref{biznorm}-(ii), for some $\sigma_1\in(0,\infty)$,
depending only on $\rho,\sigma_0,A$,
\begin{align*}
\sup_{n\geq 0} \frac{\sigma_1^nm_{2n}(0)}{(n!)^\alpha} \leq 1.
\end{align*}
We thus may apply Lemma \ref{propa1}: there is $\sigma_2\in(0,\sigma_1]$,
depending only on $\gamma,\nu,\kappa_1,\kappa_2,\sigma_1$, such that
$$
\sup_{t\geq 0} \sum_{n=0}^\infty \frac{\sigma_2^nm_{2n}(t) }{ (n!)^\alpha} \leq 3.
$$
We deduce from Lemma \ref{biznorm}-(i) that, setting $\sigma_3=\sigma_2^{1/\alpha}/2$,
$$
\sup_{t\geq 0} \intrd \exp[\sigma_3 |v|^{\rho}] f_t(\dd v) 
=
\sup_{t\geq 0} \intrd \exp[\sigma_2^{1/\alpha} |v|^{2/\alpha}/2] f_t(\dd v) 
\leq 2.3^{1/\alpha}\leq 6
$$
as desired.
\end{proof}

\bla

\end{document}